\documentclass{amsart}
\usepackage[utf8]{inputenc}
\usepackage{tikz}
\usetikzlibrary{angles,quotes}
\usepackage{float}
\usepackage{extarrows}
\usepackage{graphicx}
\usepackage{epstopdf}
\usepackage{enumerate}
\usepackage{multicol}
\usepackage{color}
\usepackage{amsmath,amssymb,amscd}
\usepackage{amsthm}
\usepackage{caption}

\usepackage{pdfpages}
\usepackage{url}

\theoremstyle{plain}
\newtheorem{thm}{Theorem}[section]
\newtheorem{lem}[thm]{Lemma}
\newtheorem{cor}[thm]{Corollary}

\newtheorem*{prop*}{Proposition}

\theoremstyle{definition}

\newtheorem{ex}[thm]{Example}

\theoremstyle{remark}

\newcommand{\sa}{\sphericalangle}

\begin{document}

\title{Minimum area isosceles containers}

\author{Gergely Kiss}
\author{J\'anos Pach}
\author{G\'abor Somlai}

\begin{abstract}
We show that every minimum area isosceles triangle containing a given triangle $T$ shares a side and an angle with $T$. This proves a conjecture of Nandakumar motivated by a computational problem. We use our result to deduce that for every triangle $T$, (1) there are at most $3$ minimum area isosceles triangles that contain $T$, and (2) there exists an isosceles triangle containing $T$ whose area is smaller than $\sqrt2$ times the area of $T$. Both bounds are best possible.
\end{abstract}

%\begin{keyword}
%isosceles triangle, special container, minimal cover
%\end{keyword}

\maketitle

\section{Introduction}\label{s1}

Given two convex bodies, $T'$ and $T$, in the plane, it is not easy to decide whether there is a rigid motion that takes $T'$ into a position where it covers $T$. Suppose, for instance, that we place a $2$-dimensional convex body $T'$ in the $3$-dimensional space, and let $T$ denote the orthogonal projection of $T'$ onto the $x$-$y$ plane. The area of $T'$ is at least as large as the area of $T$, and it looks plausible that $T'$ can be moved to cover $T$. However, the proof of this fact is far from straightforward; see~\cite{DeM86, KoT90}. As Steinhaus~\cite{St64} pointed out, it is not even clear how to decide, whether a given triangle $T'$ can be brought into a position where it covers a fixed triangle $T$. The first such algorithm was found by Post \cite{P} in 1993, and it was based on the following lemma.

\begin{lem}[Post]\label{tPost}
If a triangle $T'$ can be moved to a position where it covers another triangle $T$, then one can also find a covering position of $T'$ with a side that contains one side of $T$.
\end{lem}

In many problems, the body $T'$ is not fixed, but can be chosen from a family of possible ``containers,'' and we want to find a container which is in some sense optimal.
To find a minimum area or minimum perimeter triangle, rectangle, convex $k$-gon, or
ellipse (L\"owner-John ellipse) enclosing a given set of points are
classical problems in geometry with interesting applications in
packing and covering, approximation, convexity, computational
geometry, robotics, and elsewhere \cite{BC, BoD, D, E, Fe, I, JW, Jo, OR, R1, R2}.
Finding optimal circumscribing and inscribed simplices, ellipsoids, polytopes with a fixed number of sides or vertices, etc., are fundamental questions in optimization, functional analysis, and number theory; see e.g. \cite{GLS, KY, K,  Sc, We}.
\smallskip

Motivated by a computational problem, R. Nandakumar~\cite{N} raised the following interesting special instance of the above question: {\em Determine the minimum area of an isosceles
triangle containing a given triangle $T$.} The aim of the present note is to solve this problem and to find all triangles for which the minimum is attained. We call these triangles  {\em minimum area isosceles containers} for $T$. It is easy to verify that every triangle has at least one minimum area isosceles container (see Lemma \ref{c0.0}). However, we will see that in some cases the minimum area isosceles container is not unique.

Our main objective is to prove the following statement conjectured by Nandakumar \cite{N}.

\begin{thm}\label{t0.1}
Let $T$ be a triangle and let $T'\supseteq T$ be one of its minimum area isosceles containers. Then, $T'$ and $T$ have a side in common, and their angles at one of the endpoints of this side are equal.
\end{thm}

For any two points, $A$ and $B$, let $AB$ denote the closed segment connecting them, and let $|AB|$ stand for the length of $AB$. To unify the presentation, in the sequel we fix a triangle $T$ with vertices $A, B, C,$ and side lengths $a=|BC|$, $b=|AC|$, $c=|AB|$. If two sides are of the same length, then $T$ is the unique minimum area isosceles container of itself, so there is nothing to prove. Therefore, from now on we assume without loss of generality that $a<b<c$.
\smallskip

To establish Theorem~\ref{t0.1} and to formulate our further results, we need to introduce some special isosceles triangles associated with the triangle $ABC$, each of which shares a side and an angle with $ABC$.
\smallskip

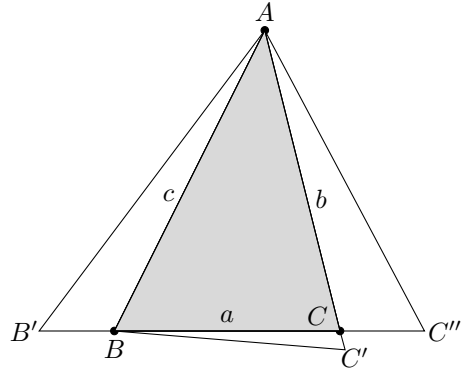
\begin{figure}[!ht]
\begin{center}
\begin{tikzpicture}
\draw (0,0) node[anchor=north]{$B$}
  -- (3,0)
  -- (2,4) node[anchor=south]{$A$}
  -- cycle;

 \draw[fill=black] (0,0) circle (0.05);
 \draw[fill=black] (3,0) circle (0.05);
 \draw[fill=black] (2,4) circle (0.05);
 \draw[fill=gray!30] (2,4)-- (0,0) --(3,0)--cycle;
\draw (1.5,0) node[anchor=south]{$a$};
\draw (2.75,2) node[anchor=north]{$b$};
\draw (0.72,2) node[anchor=north]{$c$};
\draw (0,0)
  -- (4.123,0)
  -- (2,4)
  -- cycle;
\node at (2.7,0.2 ) {$C$}; \node at (4.4,0 ) {$C''$}; \node at
(3.2,-0.35 ) {$C'$}; \draw (3,0)--(3.0675,-0.25)--(0,0);
\draw(3,0)
  -- (-1,0)
  -- (2,4)
-- cycle;

\node at (-1.2,0 ) {$B'$};

\end{tikzpicture}
\end{center}
    \caption{Special containers of the first kind $AB'C, ABC'$, and $ABC''$.}
    \label{fig2}
\end{figure}

{\bf Special containers of the first kind.} Let $B'$ denote the point on the ray $\vec{CB}$, for which $|B'C|=|AC|=b$ (see Fig.~\ref{fig2}). Analogously, let $C'$ (and $C''$) denote the
points on $\vec{AC}$ (resp., $\vec{BC}$) such that $|AC'|=c$ (resp., $|BC''|=c$). Obviously, the triangles $AB'C$, $ABC'$, and $ABC''$ are isosceles. We call them special containers of the first kind associated with $ABC$.

\smallskip

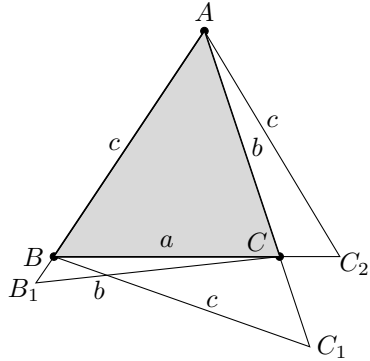
\begin{figure}[!ht]
\begin{center}
\begin{tikzpicture}
\draw[thick] (0,0) node[anchor=east]{$B$}
  -- (3,0)
  -- (2,3) node[anchor=south]{$A$}
  -- cycle;

   \draw[fill=black] (0,0) circle (0.05);
 \draw[fill=black] (3,0) circle (0.05);
 \draw[fill=black] (2,3) circle (0.05);

\draw[fill=gray!30] (2,3)-- (0,0) --(3,0)--cycle;

\draw (1.5,0) node[anchor=south]{$a$}; \draw (2.7,1.7)
node[anchor=north]{$b$}; \draw (0.8,1.7) node[anchor=north]{$c$};
\draw (2.1,-0.4) node[anchor=north]{$c$}; \draw (0.6,-0.2)
node[anchor=north]{$b$}; \draw (2.9,2) node[anchor=north]{$c$};

\draw (0,0)
  -- (3.8,0)
  -- (2,3)
  -- cycle;

\node at (2.7,0.2 ) {$C$}; \node at (4,-0.1 ) {$C_2$}; \node at
(3.7,-1.2 ) {$C_1$}; \draw (3,0)--(3.4,-1.2)--(0,0);
\draw(0,0)
  -- (-0.24,-0.35)
  -- (3,0)
-- cycle;

\node at (-0.4,-0.45 ) {$B_1$};
\end{tikzpicture}
\end{center}
    \caption{Special containers of the second kind $AB_1C$, $ABC_1$, and $ABC_2$.}
    \label{fig2.1}
\end{figure}

{\bf Special containers of the second kind.} Let $B_1$ denote the point on the ray $\vec{AB}$, different from $A$, for which $|B_1C|=|AC|=b$ (see Figure \ref{fig2.1}). Analogously, let $C_1$ (resp., $C_2$) denote the point on $\vec{AC}$ (resp., $\vec{BC}$) for which $|BC_1|=|AB|=c$ and $C_1\neq A$ (resp., $|AC_2|=|AB|=c$ and $C_2\neq B$). The triangles $AB_1C$, $ABC_1$, and $ABC_2$ are called the special containers of the second kind associated with $ABC$.

\smallskip

\begin{figure}[!ht]
\begin{center}
\begin{tikzpicture}
\draw (0,0) node[anchor=north]{$B$}
  -- (3,0) node[anchor=north]{$C$}
  -- (2,3) node[anchor=west]{$A$}
  -- cycle;

 \draw[fill=black] (0,0) circle (0.05);
 \draw[fill=black] (3,0) circle (0.05);
 \draw[fill=black] (2,3) circle (0.05);
 \draw[fill=gray!30] (2,3)-- (0,0) --(3,0)--cycle;

\draw (1.5,0)--(1.5, 4.5)--(0,0); \draw (1.5,4.5)-- (2,3); \draw (1,
1.5)--(3.25,0); \draw (2.5, 1.5)--(-2,0);
\draw (0,0)
  -- (3.25,0)
  -- (2,3)
  -- cycle;

\node at (3.5,0 ) {$\overline{C}$}; \node at (1.5,4.8 )
{$\overline{A}$};
\draw(3,0)
  -- (-2,0)
  -- (2,3)
-- cycle;

\node at (-2.3,0 ) {$\overline{B}$};

\draw[thick] (5.5,0) node[anchor=north]{$B$}
  -- (7,3)node[anchor=west]{$A$}
  -- (6.5,0) node[anchor=north]{$C$}
  -- cycle;

\draw[fill=black] (5.5,0) circle (0.05);
\draw[fill=black] (7,3) circle (0.05);
\draw[fill=black] (6.5,0) circle (0.05);
\draw[fill=gray!30] (5.5,0)-- (7, 3) --(6.5,0)--cycle;
\draw (6.5,0)--(9.25, 0)--(7,3);
\draw (9.25, 0)-- (6.25,1.5);

\node at (9.45,0) {$\overline{C}$};
\end{tikzpicture}
\end{center}
    \caption{Special containers of the third kind in the acute and in the non-acute cases.}
    \label{fig2.2}
\end{figure}
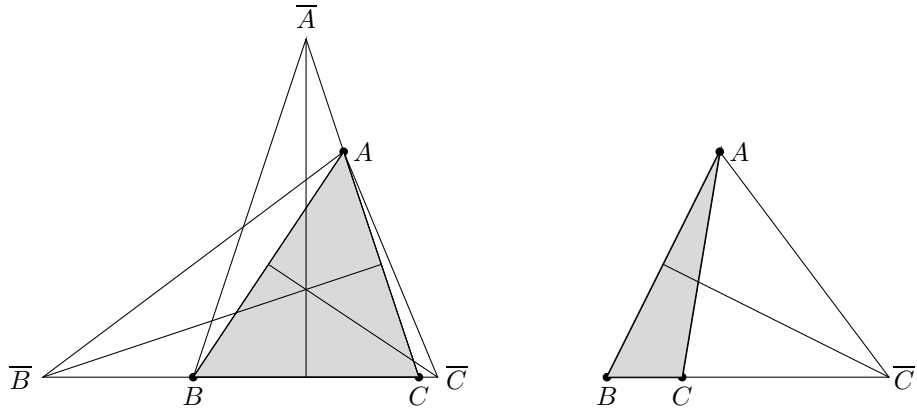

{\bf Special containers of the third kind.} Let $\overline{A}$ be the intersection of the perpendicular bisector of $BC$ and the line $AC$. Since we have $b=|AC|<|AB|=c$, the point
$\overline{A}$ lies outside of $ABC$. Analogously, denote by $\overline{B}$ (resp., $\overline{C}$) the intersection of the perpendicular bisector of $AC$ (resp. $AB$) and the line $BC$. Note that $\overline{A}BC$ and $A\overline{B}C$ do not contain $ABC$ if $\sa BCA\ge 90  ^{\circ}$. The triangles $\overline{A}BC$, $A\overline{B}C$, and $AB\overline{C}$ are called special
containers of the third kind associated with $ABC$, provided that they contain $ABC$.
Thus, if $ABC$ is acute, then it has three special containers of the third kind. Otherwise, it has only one (see Figure \ref{fig2.2}).

\smallskip

All special containers share a common angle and a common side with the original triangle $ABC$. Obviously, there is no other isosceles container having the same property. Indeed, for each vertex of $ABC$, there are at most 3 isosceles triangles that share this vertex and the angle at this vertex with $ABC$, and also have a common side with $ABC$.

Therefore, Theorem \ref{t0.1} is an immediate corollary of the following statement.

\begin{thm}\label{t0.2}
All minimum area isosceles containers for a triangle are special containers of the first kind, or of the second kind, or of the third kind.
\end{thm}

Whenever a minimum area isosceles container of a triangle is acute, we can be more specific.

\begin{thm}\label{prop1.3}
If a minimum area isosceles container of a triangle is acute, then it is a special container of the first kind.
\end{thm}

One is tempted to believe that if a triangle is acute, then all of its minimum area isosceles containers are acute and, hence, all of them are special containers of the first kind. However, this is not the case: Example~\ref{ex1} demonstrates that there are acute triangles with obtuse minimum area isosceles containers. As all special containers of the first kind and the third kind of an acute triangle are acute, Theorem~\ref{prop1.3} implies the following statement.

\begin{cor}\label{corp2}
A minimum area isosceles container for an acute triangle is obtuse if and only if it is a special container of the second kind.
\end{cor}

It follows from Theorem \ref{t0.2} that every triangle has at most $9$ minimum area isosceles containers: at most $3$ special containers of each kind. In the next section, we prove that there are no minimum area isosceles triangles of the third kind (see Lemma~\ref{p2}). Thus, every triangle can have at most $6$ minimum area isosceles triangles. In fact, this bound can be further reduced to $3$.

\begin{thm}\label{peq3}
Every non-isosceles triangle $ABC$ has at most $3$ minimum area isosceles containers, $AB'C$, $ABC'$, and $AB_1C$. In particular, every minimum area isosceles container is a special container of the first or the second kind.

There is a unique triangle $T^*$, up to similarity, which has precisely $3$ different minimum area isosceles containers. Its angles are $\alpha^*\approx 41.831452 ^{\circ}, 2\alpha^*,$ and $180 ^{\circ}-3\alpha^*$, where $\alpha^*$ is the unique solution of  $\sin(\alpha)\sin (2\alpha)-\sin^2(3\alpha)=0$ in the interval $[36 ^{\circ}, 45^{\circ}]$.
\end{thm}

Finally, we discuss how large the area of a minimum area isosceles container for a triangle $T$ can be relative to the area of $T$. We also consider the same question for special containers of the first kind.

\begin{thm}\label{thmrarea}
\begin{enumerate}[(a)]
\item\label{item16a} Every triangle of area $1$ has an isosceles container whose area is smaller than $\sqrt{2}$.
\item\label{item16b} Every triangle of area $1$ has a special container of the first kind, whose area is smaller than $\frac{1+\sqrt{5}}{2}$.
\end{enumerate}
\smallskip

Both bounds are best possible.
\end{thm}

As (b) is best possible, there exists a triangle of area 1 for which every special container of the first kind has area larger than $\frac{1+\sqrt{5}}{2}$. Therefore, by (a), none of its special containers of the first kind can be a minimum area isosceles container. This disproves an earlier conjecture of Nandakumar, according to which every triangle $T$ admits a minimum area isosceles container which is a special container of the first kind.
\smallskip

Our paper is organized as follows. In Section \ref{s2}, we prove some useful inequalities for the areas of special containers. In particular, we prove Theorem \ref{thmrarea} (b). In Section \ref{s3}, we establish some elementary properties of minimum area isosceles containers. Section \ref{smain} contains the proofs of Theorems \ref{t0.2}, and \ref{prop1.3}. Finally, Theorem \ref{peq3}, and \ref{thmrarea} (a) are proved in Section \ref{uj}.

\section{Preliminaries---Proof of Theorem~\ref{thmrarea} (\rm{b})}\label{s2}
In this section, we collect some basic facts about special containers and establish Theorem~\ref{thmrarea} (b).

First, we consider special containers of the first kind, because their areas can be easily compared to the area of the triangle $ABC$. As everywhere else, we assume that the side lengths of $ABC$ satisfy $a<b<c$. The area of $ABC$ is denoted by $t(ABC)$.

\begin{lem}\label{p0.1} For any non-isosceles triangle $ABC$,  we have
\begin{enumerate}
    \item[(a)] $t(ABC'')>t(ABC')$,
    \item[(b)] $t(AB'C)>t(ABC')$ (resp., $t(AB'C)\geq t(ABC')$) \\ if and only if $b^2>ac$ (resp., $b^2\geq ac$ ).
\end{enumerate}
\end{lem}

%\begin{proof}
\noindent{\bf Proof.} Let $m$ be the length of the altitude of $ABC$ perpendicular
to the side $BC$. We have $t(ABC)=\frac{a\cdot m}{2}$ and
$t(AB'C)=\frac{b\cdot m}{2}$. Thus, the ratio $t(AB'C)/t(ABC)=
b/a$.

\noindent Similar arguments show that
$$\frac{t(ABC')}{t(ABC)}= \frac{c}{b} \ \ \textrm{ and } \ \ \frac{t(ABC'')}{t(ABC)}=\frac{c}{a}.$$

\begin{enumerate}
    \item[{\it(a)}]
Since $b>a$, we have $t(ABC'')>t(ABC')$.
   \item[{\it (b)}]
    Straightforward.  \hfill $\Box$
\end{enumerate}

\medskip
Next, using Lemma \ref{p0.1}, we determine the supremum of the ratio of the area of the smallest isosceles containers of the first kind to the area of the original triangle $ABC$. This will also provide an upper bound for the ratio of the area of a smallest area {\em isosceles} containers to the area of the original triangle. This fact will be used in the sequel.

\medskip

\noindent{\bf Proof of Theorem \ref{thmrarea} (b).}
Let $r^*_1$ denote the supremum of the ratio of the area of a smallest container of the first kind associated with $ABC$ and the area of $ABC$, over all triangles $ABC$ with the above property.
We show that $r^*_1= \frac{1+\sqrt{5}}{2}$. Suppose without loss of generality that $a=1$. By our assumptions and the triangle inequality, we have $1<b<c<b+1$.
Let
$$r(b,c):=\begin{cases}
b& \textrm{if } b^2\le c,\\
c/b & \textrm{if } b^2> c.
\end{cases}$$
By Lemma \ref{p0.1}, $r(b,c)=\min\Big(
\frac{t(AB'C)}{t(ABC)},\frac{t(ABC')}{t(ABC)}\Big)$, and
\begin{equation*}\label{eq_fa} r_1^*=\sup_{1<b<c<b+1}
r(b,c).%=\frac{1+\sqrt{5}}{2}.
\end{equation*}
If
$b^2\le c<b+1$, then $r(b,c)=b < \frac{1+\sqrt{5}}{2}$.

If $b^2>c$ and $b < \frac{1+\sqrt{5}}{2}$, then
$r(b,c)=\frac{c}{b}<b < \frac{1+\sqrt{5}}{2}$. Otherwise,
if
$b\ge\frac{1+\sqrt{5}}{2}$ (and, hence, $b^2>c$), then
$r(b,c)=\frac{c}{b}<\frac{b+1}{b}=1+\frac{1}{b}\le 1+\frac{2}{1+\sqrt{5}}=\frac{1+\sqrt{5}}{2}$.
Thus, we obtain that $r(b,c)<r_1^*\le \frac{1+\sqrt{5}}{2}$ , for all $1<b<c<b+1$, i.e, the supremum $r_1^*$ is not attained for any triangle $ABC$.

The supremum of $r(b,c)$, restricted to the parabola arc $c=b^2<b+1$ in the $(b,c)$ plane, is $\frac{1+\sqrt{5}}{2}$. Since every point $(b,c)$ of this arc corresponds to a triangle with side lengths $1,b,c$, we obtain that $r^*_1=\frac{1+\sqrt{5}}{2}$, as required.
\qed

\bigskip

We show that for any special container of the third kind, there
exists a special container of the second kind whose area is smaller.

\begin{lem}\label{p2}
For any triangle $ABC$ we have
$ t(AB_1C) < t(AB\overline{C})$. If $ABC$ is acute, then we also have $t(ABC_1) < t(A\overline{B}C)$ and $t(ABC_2) < t(\overline{A}BC)$.
Thus, a minimum area isosceles container can never be a special container of the third kind.
\end{lem}

\noindent{\bf Proof.} We verify only the first inequality; the other two statements can be shown analogously.

Assign planar coordinates to the points.
We can assume without loss of generality that $A=(0,0)$, $B=(2,0)$,  $C=(p,q)$, and $\overline{C}=(1,d)$. Then
$t(AB\overline{C})=d$.
The equation of the line passing through $B$, $C$, and
$\overline{C}$ is $dx+y=2d$. Taking $p=1+s$ for some $0<s<1$, we have $q=d(1-s)$ and $B_1=(2(1+s),0)$. Hence,
$t(AB_1C)=d(1-s^2)<d=t(AB\overline{C})$.
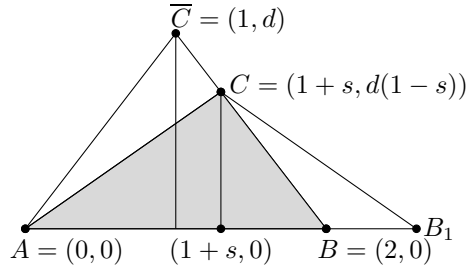
\begin{figure}[!ht]
\begin{center}
\begin{tikzpicture}

\draw (0,0)--(2,2.6)--(4,0);

 \draw[fill=black] (4,0) circle (0.05);
 \draw[fill=black] (0,0) circle (0.05);
 \draw[fill=black] (2.6,1.82) circle (0.05);
 \draw[fill=gray!30] (0,0)-- (2.6,1.82) --(4,0)--cycle;

\draw (0,0)--(2.6,1.82)--(5.2,0)--cycle; \draw (2,2.6)--(2,0); \draw
(2.6,1.82)--(2.6,0); \node at (0.55,-0.3) {$A=(0,0)$}; \node at
(4.65,-0.3) {$B=(2,0)$}; \node at (2.6,-0.3) {$(1+s,0)$}; \node at
(2.7,2.8) {$\overline{C}=(1,d)$};

\node at (5.5,0) {$B_1$}; \node at (4.3,1.9) {$C=(1+s,d(1-s))$};
\draw[fill=black] (0,0) circle (0.05); \draw[fill=black] (2.6,0)
circle (0.05); \draw[fill=black] (5.2,0) circle (0.05);
\draw[fill=black] (4,0) circle (0.05); \draw[fill=black] (2.6,1.82)
circle (0.05); \draw[fill=black] (2,2.6) circle (0.05);
\end{tikzpicture}
\caption{Proof of Lemma \ref{p2}.}
    \label{fig4.2}
\end{center}
\end{figure}
\hfill $\Box$

\section{Four useful lemmas}\label{s3}

The aim of this section is to prepare the ground for the proofs of the main results that will be given in the next two sections.

First, we show that the problem is well-defined, that is, for every triangle $ABC$, there is at least one isosceles triangle containing $ABC$, whose area is smaller than or equal to the area of any other isosceles container.
\begin{lem}\label{c0.0}
Every triangle $ABC$ has at least one minimum area isosceles container.
 \end{lem}

\noindent{\bf Proof.}
It follows from Theorem \ref{thmrarea} (b) that the area of a minimum area isosceles container is at most $\frac{1+\sqrt{5}}{2}$ times larger than the area of the original triangle $ABC$. Therefore, the vertices of any minimum area isosceles container must lie within a bounded distance from
$ABC$, and the statement follows by a standard compactness argument.  \hfill  $\Box$

\begin{lem}\label{c0}
Let $ABC$ be a triangle and $SPR$ a minimum area
isosceles container for $ABC$. Then $A,B,C$ must lie on the boundary of $SPR$, and each side of $SPR$ contains at least one of them.
\end{lem}

\noindent{\bf Proof.}
First, we show that each side of $SPR$ contains a vertex of $ABC$. Indeed, if one of the sides did not contain any vertex of $ABC$, then we could slightly move it, parallel to itself, to obtain a smaller isosceles container.

Assume next, for contradiction, that $A$, say, does not lie on the boundary of $SPR$. Then we could slightly rotate $ABC$ about $B$ or $C$, to bring it into a position where only one of its vertices lies on the boundary of $SPR$. However, in that case, at least one of the sides of $SPR$ would contain no vertex of $ABC$.
\hfill   $\Box$
\medskip

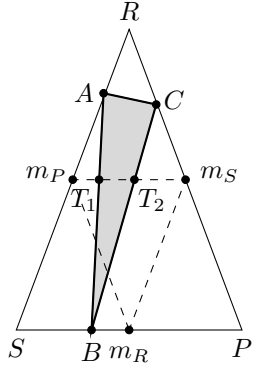
\begin{figure}[H]
\begin{center}
\begin{tikzpicture}
\draw (0,0) node[anchor=north]{$S$}
  -- (3,0) node[anchor=north]{$P$}
  -- (1.5,4) node[anchor=south]{$R$}
  -- cycle;

 \draw[fill=gray!30] (1.16,3.15)-- (1,0) --(1.86,3)--cycle;

 \draw[fill=black] (1.5,0) circle (0.05);
 \draw[fill=black] (2.25,2) circle (0.05);
 \draw[fill=black] (0.75,2) circle (0.05);

\draw[fill=black] (1.57,2) circle (0.05); \draw (1.8,2)
node[anchor=north]{$T_2$}; \draw[fill=black] (1.1,2) circle (0.05);
\draw (0.9,2) node[anchor=north]{$T_1$};

\draw (1.5,-0.1) node[anchor=north]{$m_{R}$}; \draw (2.7,2.3)
node[anchor=north]{$m_{S}$}; \draw (0.4,2.3)
node[anchor=north]{$m_{P}$};
%\draw[dashed] (1.5, 0)--(0,0)--(0.75,2);
\draw[dashed] (1.5, 0) -- (0.75,2) --(2.25,2)--cycle;
%\draw[dashed]  (0.75,2)-- (1.5,4) --(2.25,2);

\draw[fill=black] (1,0) circle (0.05);
 \draw[fill=black] (1.16,3.15) circle (0.05);
 \draw[fill=black] (1.86,3) circle (0.05);

 \draw (1,-0.05) node[anchor=north]{$B$};
\draw (0.9,3.4) node[anchor=north]{$A$}; \draw (2.1,3.3)
node[anchor=north]{$C$};

\draw[thick] (1,0)-- (1.16,3.15) -- (1.86,3)--cycle;

\end{tikzpicture}
\end{center}
 \caption{Illustration for the proof of Lemma \ref{clm2}.}
    \label{fig7}
\end{figure}

Let $SPR$ be an isosceles triangle and let $m_{S}, m_{P}, m_{R}$ denote the midpoints of the sides $PR$, $SR$, and $SP$, respectively.
The boundary of $SPR$ splits into three polygonal pieces,
$\widehat{m_Pm_R}, \widehat{m_Rm_S}$ and $\widehat{m_Sm_P}$, each of
which consists of two closed line segments.
Namely,
$$ \widehat{m_Pm_R}= m_PS ~ \cup ~ Sm_R,$$ $$\widehat{m_Rm_S}=m_RP ~\cup ~ P m_S,$$  $$\widehat{m_Sm_P}= m_SR ~\cup~ Rm_P.$$
See Figure \ref{fig7}.

\begin{lem}\label{clm2}
Let $ABC$ be a triangle and $SPR$ a minimum area
isosceles container for $ABC$. Then, each of the closed polygonal pieces
$\widehat{m_Pm_R}, \widehat{m_Rm_S}$, and $\widehat{m_Sm_P}$ contains
precisely one vertex of $ABC$.
\end{lem}
\smallskip

\noindent{\bf Proof.} By Lemma \ref{c0}, the vertices $A,B,C$ lie on the boundary of $SPR$. Suppose for contradiction that the closed polygonal piece $\widehat{m_Pm_S}$ contains two
vertices of $ABC$. We may and do assume without loss of generality that these vertices are $A$ and $C$.

Let $T_1$ and $T_2$ denote the intersection points of the segment $m_Pm_S$
with $AB$ and $CB$, respectively.
The quadrilateral $CAT_1T_2\subseteq Rm_Pm_S$, so that $t(CAT_1T_2)\le t(Rm_Pm_S)$. Since $|T_1T_2|\le|m_Pm_S|$, we have $t(T_2 T_1B)\le t(m_Sm_Pm_R)$. Consequently, we get $$t(ABC)\le t(Rm_Pm_S) + t(m_Sm_Pm_R)=\frac{1}{2}t(SPR).$$
Equality holds if and only if $A=R,B=S, C=m_S$ or
$A=m_P,B=P, C=R$.

On the other hand,
by Theorem \ref{thmrarea} (b), we obtain $t(SPR)<\frac{1+\sqrt{5}}{2}t(ABC)$, the desired contradiction. \hfill $\Box$

\begin{lem}\label{lcommon} Let $ABC$ be a triangle and $SPR$ a minimal area isosceles container for $ABC$. Then
$ABC$ and $SPR$ have a common vertex.
\end{lem}

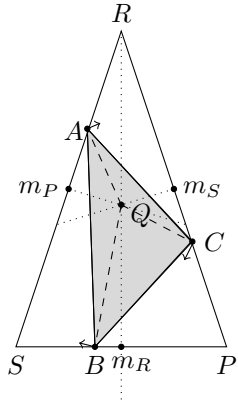
\begin{figure}[!ht]
\begin{center}
\begin{tikzpicture}[scale=0.7]

\draw (0,0) node[anchor=north]{$S$}
  -- (4,0) node[anchor=north]{$P$}
  -- (2,6) node[anchor=south]{$R$}
  -- cycle;

 \draw[fill=gray!30] (1.355,4.145)-- (1.5,0) --(3.35,2)--cycle;

\node[anchor=north] at (2.2,0) {$m_R$}; \node[anchor=east] at (1,3)
{$m_P$}; \node[anchor=west] at (3,3) {$m_S$}; \draw[fill=black]
(2,0) circle (0.05); \draw[fill=black] (1,3) circle (0.05);
\draw[fill=black] (3,3) circle (0.05);

\draw (1.355,4.145)--(1.5,0)--(3.33,2)--cycle;
 \draw[fill=black] (1.355,4.145) circle (0.05);
 \draw[fill=black] (1.5,0) circle (0.05);
 \draw[fill=black] (3.35,2) circle (0.05);

\draw (1.1,4.45 ) node[anchor=north]{$A$}; \draw (1.5, 0 )
node[anchor=north]{$B$}; \draw (3.4,2) node[anchor=west]{$C$};

\draw[->] (1.355,4.145) --(1.6, 4.25);
\draw[->] (1.5,0 )--(1.2, 0.075);
\draw[->] (3.35,1.95 )--(3.2, 1.65);

\draw[fill=black] (2, 2.7) circle
(0.05);

\draw[dotted] (2, -1)--(2, 6);

\draw[dotted] (3, 3)--(0.7, 2.3);
\draw[dotted] (1, 3)--(3.3, 2.3);

\draw  (2, 2.5)  node[anchor=west] {$Q$};
\draw[dashed] (2,2.7)--(1.355,4.145); \draw[dashed]
(2,2.7)--(1.5,0); \draw[dashed] (2,2.7)--(3.33,2);

\end{tikzpicture}
\end{center}
 \caption{Illustration for the proof of Lemma \ref{lcommon}.}
    \label{figcase1}
\end{figure}

\noindent{\bf Proof.}
By Lemma \ref{c0}, the points $A,B,C$ must lie on the boundary $SPR$. Suppose for contradiction that none of them is a vertex of $SPR$. In view of Lemmas \ref{c0} and \ref{clm2}, there are two possibilities: each of the segments $m_PR,m_RS, m_SP$ contains precisely one vertex of $ABC$, or each of the segments $m_PS,m_RP, m_SR$ contains precisely one vertex of $ABC$. Suppose without loss of generality that $A\in m_PR,B\in m_RS, C\in m_SP$, as in Figure \ref{figcase1}.

Let $Q$ denote the center of the circle circumscribed around $SPR$. It is easy to see that a small clockwise rotation about $Q$ will take $ABC$ into a position such that all of its vertices lie in the interior of the triangle $SPR$. This contradicts the minimality of $SPR$ and Lemma \ref{c0}.
\hfill $\Box$

\section{Nandakumar's conjecture---Proofs of Theorems~\ref{t0.2} and~\ref{prop1.3}}\label{smain}

We start with the proof of Theorem \ref{t0.2}, which immediately implies Nandakumar's conjecture (Theorem~\ref{t0.1}).
\medskip

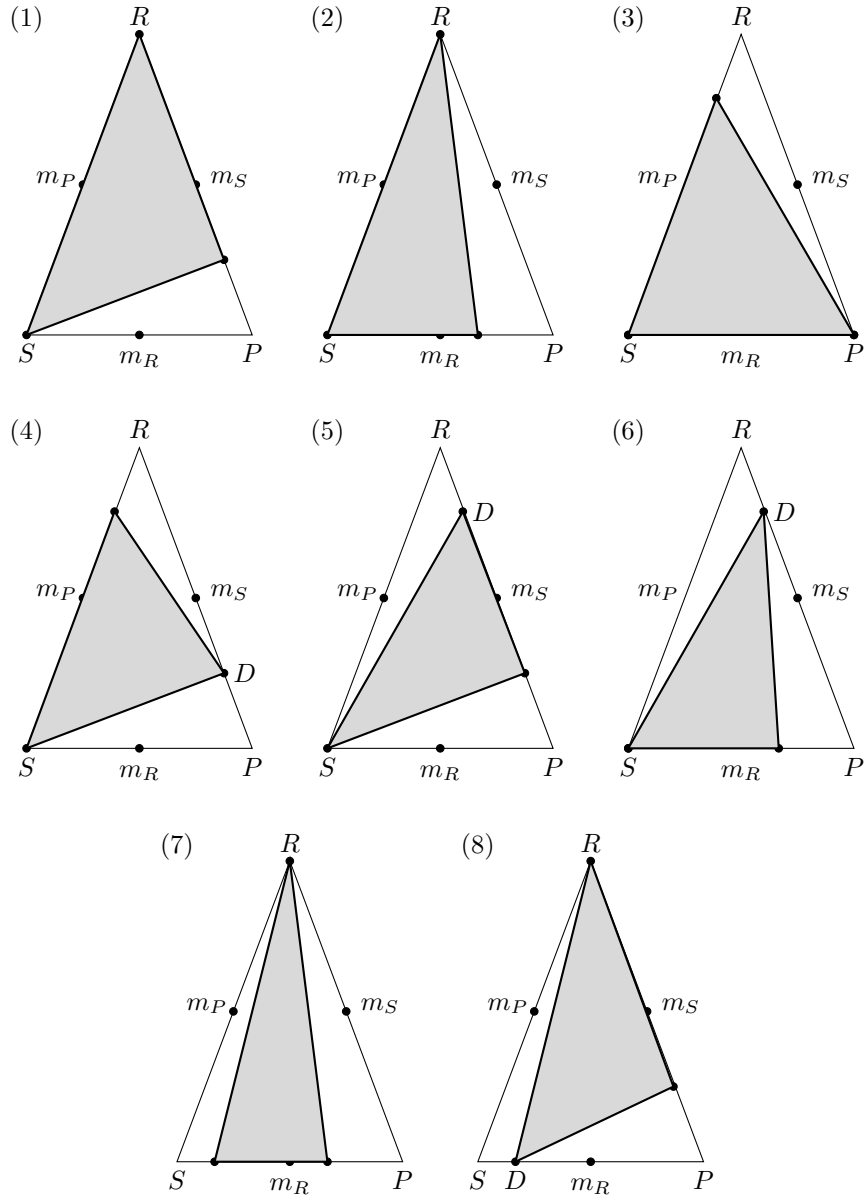
\begin{figure}[!ht]
\begin{center}
\begin{tikzpicture}
%(1)
\draw (0,0) node[anchor=north]{$S$}
  -- (3,0) node[anchor=north]{$P$}
  -- (1.5,4) node[anchor=south]{$R$}
  -- cycle;
\draw (0,3.9) node[anchor=south]{$(1)$};
 \draw[fill=black] (1.5,0) circle (0.05);
 \draw[fill=black] (2.25,2) circle (0.05);
 \draw[fill=black] (0.75,2) circle (0.05);

\draw (1.5,-0.1) node[anchor=north]{$m_{R}$}; \draw (2.7,2.3)
node[anchor=north]{$m_{S}$}; \draw (0.4,2.3)
node[anchor=north]{$m_{P}$};

\draw[fill=black] (0,0) circle (0.05);
 \draw[fill=black] (1.5,4) circle (0.05);
 \draw[fill=black] (2.625,1) circle (0.05);

\draw[thick, fill=gray!30] (0,0)-- (1.5,4) -- (2.625,1)-- cycle;

%(2)
\draw (4,0) node[anchor=north]{$S$}
  -- (7,0) node[anchor=north]{$P$}
  -- (5.5,4) node[anchor=south]{$R$}
  -- cycle;

\draw (4,3.9) node[anchor=south]{$(2)$};
 \draw[fill=black] (5.5,0) circle (0.05);
 \draw[fill=black] (6.25,2) circle (0.05);
 \draw[fill=black] (4.75,2) circle (0.05);

\draw (5.5,-0.1) node[anchor=north]{$m_{R}$}; \draw (6.7,2.3)
node[anchor=north]{$m_{S}$}; \draw (4.4,2.3)
node[anchor=north]{$m_{P}$};

\draw[fill=black] (4,0) circle (0.05);
 \draw[fill=black] (5.5,4) circle (0.05);
 \draw[fill=black] (6,0) circle (0.05);

\draw[thick, fill=gray!30] (4,0)-- (5.5,4) -- (6,0)-- cycle;

%(3)
\draw (8,0) node[anchor=north]{$S$}
  -- (11,0) node[anchor=north]{$P$}
  -- (9.5,4) node[anchor=south]{$R$}
  -- cycle;
 \draw (8,3.9) node[anchor=south]{$(3)$};
 \draw[fill=black] (8,0) circle (0.05);
 \draw[fill=black] (10.25,2) circle (0.05);
 \draw[fill=black] (11,0) circle (0.05);

\draw (9.5,-0.1) node[anchor=north]{$m_{R}$}; \draw (10.7,2.3)
node[anchor=north]{$m_{S}$}; \draw (8.4,2.3)
node[anchor=north]{$m_{P}$};

\draw[fill=black] (8,0) circle (0.05);
 \draw[fill=black] (9.17,3.15) circle (0.05);
 \draw[fill=black] (11,0) circle (0.05);

\draw[thick, fill=gray!30] (8,0)-- (9.17,3.15) -- (11,0)-- cycle;

%%%%%%%%%%
%(5)
\draw (0,-5.5) node[anchor=north]{$S$}
  -- (3,-5.5) node[anchor=north]{$P$}
  -- (1.5,-1.5) node[anchor=south]{$R$}
  -- cycle;
\draw (0,-1.6) node[anchor=south]{$(4)$};

 \draw[fill=black] (1.5,-5.5) circle (0.05);
 \draw[fill=black] (2.25,-3.5) circle (0.05);
 \draw[fill=black] (0.75,-3.5) circle (0.05);

\draw (1.5,-5.6) node[anchor=north]{$m_{R}$}; \draw (2.7,-3.2)
node[anchor=north]{$m_{S}$}; \draw (0.4,-3.2)
node[anchor=north]{$m_{P}$};

\draw[fill=black] (0,-5.5) circle (0.05);
 \draw[fill=black] (1.17,-2.35) circle (0.05);
 \draw[fill=black] (2.625,-4.5) circle (0.05);

\draw[thick, fill=gray!30] (0,-5.5)-- (1.17,-2.35) -- (2.625,-4.5)--
cycle;

%(6)
\draw (4,-5.5) node[anchor=north]{$S$}
  -- (7,-5.5) node[anchor=north]{$P$}
  -- (5.5,-1.5) node[anchor=south]{$R$}
  -- cycle;

\draw (4,-1.6) node[anchor=south]{$(5)$};
 \draw[fill=black] (5.5,-5.5) circle (0.05);
 \draw[fill=black] (6.25,-3.5) circle (0.05);
 \draw[fill=black] (4.75,-3.5) circle (0.05);

\draw (5.5,-5.6) node[anchor=north]{$m_{R}$}; \draw (6.7,-3.2)
node[anchor=north]{$m_{S}$}; \draw (4.4,-3.2)
node[anchor=north]{$m_{P}$}; \draw (2.625,-4.5)
node[anchor=west]{$D$};

\draw[fill=black] (6.625,-4.5) circle (0.05);
 \draw[fill=black] (5.8,-2.35) circle (0.05);
 \draw[fill=black] (4,-5.5) circle (0.05);

 \draw (5.8,-2.35) node[anchor=west]{$D$};

\draw[thick, fill=gray!30] (4,-5.5)-- (5.8,-2.35) -- (6.625,-4.5)--
cycle;

%(7)
\draw (8,-5.5) node[anchor=north]{$S$}
  -- (11,-5.5) node[anchor=north]{$P$}
  -- (9.5,-1.5) node[anchor=south]{$R$}
  -- cycle;
 \draw (8,-1.6) node[anchor=south]{$(6)$};
 \draw[fill=black] (8,-5.5) circle (0.05);
 \draw[fill=black] (10.25,-3.5) circle (0.05);

\draw[fill=black] (9.5,-5.6) node[anchor=north]{$m_{R}$}; \draw
(10.7,-3.2) node[anchor=north]{$m_{S}$}; \draw (8.4,-3.2)
node[anchor=north]{$m_{P}$}; \draw[fill=black] (8,-5.5) circle
(0.05);
 \draw[fill=black] (9.8,-2.35) circle (0.05);
 \draw[fill=black] (10,-5.5) circle (0.05);

 \draw (9.8,-2.35) node[anchor=west]{$D$};

\draw[thick, fill=gray!30] (8,-5.5)-- (10,-5.5) -- (9.8,-2.35)--
cycle;

%(8)
\draw (2,-11) node[anchor=north]{$S$}
  -- (5,-11) node[anchor=north]{$P$}
  -- (3.5,-7) node[anchor=south]{$R$}
  -- cycle;
\draw (2,-7.1) node[anchor=south]{$(7)$};
 \draw[fill=black] (3.5,-11) circle (0.05);
 \draw[fill=black] (4.25,-9) circle (0.05);
 \draw[fill=black] (2.75,-9) circle (0.05);

\draw (3.5,-11.1) node[anchor=north]{$m_{R}$}; \draw (4.7,-8.7)
node[anchor=north]{$m_{S}$}; \draw (2.4,-8.7)
node[anchor=north]{$m_{P}$};

\draw[fill=black] (2.5,-11) circle (0.05);
 \draw[fill=black] (3.5,-7) circle (0.05);
 \draw[fill=black] (4,-11) circle (0.05);

\draw[thick, fill=gray!30] (2.5,-11)-- (4,-11) -- (3.5,-7)-- cycle;

%%%%%%%%%%%%%%%%
%(kimaradt)
\draw (6,-11) node[anchor=north]{$S$}
  -- (9,-11) node[anchor=north]{$P$}
  -- (7.5,-7) node[anchor=south]{$R$}
  -- cycle;
\draw (6,-7.1) node[anchor=south]{$(8)$};
 \draw[fill=black] (7.5,-11) circle (0.05);
 \draw[fill=black] (8.25,-9) circle (0.05);
 \draw[fill=black] (6.75,-9) circle (0.05);

\draw (7.5,-11.1) node[anchor=north]{$m_{R}$}; \draw (8.7,-8.7)
node[anchor=north]{$m_{S}$}; \draw (6.4,-8.7)
node[anchor=north]{$m_{P}$}; \draw (6.5,-11)
node[anchor=north]{$D$};

\draw[fill=black] (6.5,-11) circle (0.05);
 \draw[fill=black] (7.5,-7) circle (0.05);
 \draw[fill=black] (8.6,-10) circle (0.05);

\draw[thick, fill=gray!30] (6.5,-11)-- (8.6,-10) -- (7.5,-7)--
cycle;

\end{tikzpicture}
\end{center}
 \caption{The 8 cases up to symmetry. Triangle $ABC$ is shaded.}

    \label{figesetek}
\end{figure}

\noindent{\bf Proof of Theorem \ref{t0.2}.} Let $SPR$ be a minimum area container for $ABC$ with apex $R$. By Lemma \ref{c0}, $A,B$, and
$C$ are on the boundary of $SPR$, and, by Lemma \ref{lcommon}, the triangles $SPR$ and $ABC$ share a
vertex. Using Lemma \ref{clm2} under the assumption that $ABC\ne
SPR$, we can distinguish 8 cases, up to symmetry (see
Figure \ref{figesetek}). Cases (1)--(3) represent those instances when
$ABC$ and $SPR$ have two common vertices. In these cases, $SPR$ is a special
container of the first, the second, and the third kind,
respectively, so we are done.

In the remaining cases, $ABC$ and $SPR$ have only one vertex in common.
In cases (4)--(6), this vertex is a base vertex (say, $S$) of $SPR$.
Finally, in cases (7)--(8), $R$ is the unique common vertex of $ABC$ and
$SPR$. It is sufficient to show that in cases (4)--(8), the area of $SPR$ is not
minimal.

First, we discuss cases (5)--(8). Case (4) is more delicate and is left to the end of the proof.

Cases (5) and (6) are analogous. Let $D$ denote the vertex of $ABC$ lying on $RP$. In
both cases, we have $ABC\subseteq SPD$. Clearly, $SPR$ is a special container
of the third kind associated with $SPD$, and by Lemma
\ref{p2}, it cannot be minimal. Since every container for $SPD$ is also a container for $ABC$, we
conclude that $SPR$ is not a minimum area container for $ABC$.

In case (7), we can find an isosceles triangle with apex $R$ which contains $ABC$ and whose base is properly contained in $PS$. Thus, $SPR$ was not minimal.

In case (8), one vertex of $ABC$ is $R$, another (denoted by $D$) belongs to $Sm_R$, and the third lies on $Pm_S$. Since $SPR$ is an
isosceles triangle, we have $\sa RDS\ge 90 ^{\circ}$. Hence, $ABC$ can be slightly
rotated about $R$ so that it remains within $SPR$, which leads to a contradiction.
\smallskip

It remains to handle case (4). We distinguish two subcases.
Denote the apex angle $\sa SRP$ by $\delta$. If $\delta\ge 60 ^{\circ}$,
then we can rotate $ABC$ about $S$. Indeed, vertex $D$ of $ABC$ belongs to $m_SP$, while the base of the altitude belonging to $PR$ lies on $m_SR$.
Hence, we have $\sa SDP\ge 90 ^{\circ}$, and the image of
$ABC$ through a small clockwise rotation %with small negative angle
about $S$ is still contained in $SPR$. Therefore, in this case, $SPR$ cannot be minimal either.

Therefore, from now on we assume $\delta< 60  ^{\circ}$.
Choose a suitable coordinate system, in which the vertices of $ABC$ are
$(0,0)$, $(s,0)$, and $D=(p,q)$. We also have $S=(0,0)$ and $R=(s+x,0)$ for some $x>
0$. Since $\delta<60 ^{\circ}<90 ^{\circ}$, vertex $D$ is to the left of $R$, that is, $p<s+x$.

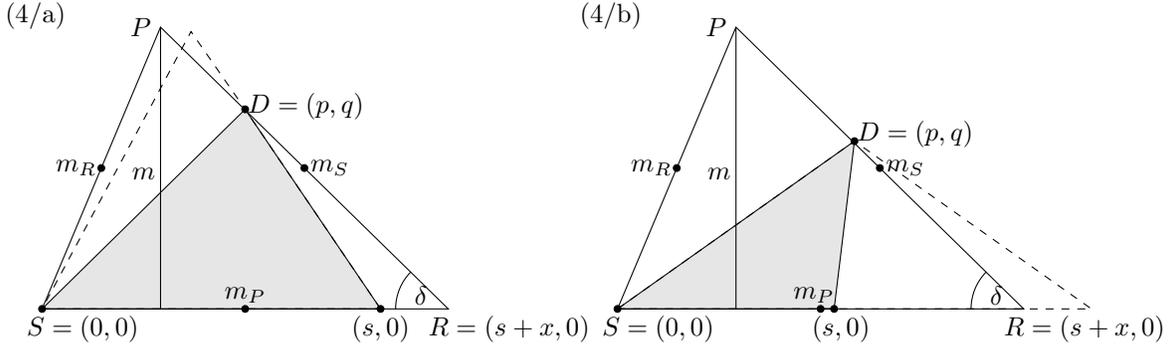
\begin{figure}[!ht]
\begin{center}
\begin{tikzpicture}[scale=0.90]
\draw (-0.1,4) node[anchor=south]{(4/a)}; \draw (8.4,4)
node[anchor=south]{(4/b)}; \draw[fill=gray!20]
(0,0)--(5,0)--(3,2.95)--cycle; \draw[dashed]
(5,0)--(3,2.95)--(2.2,4.1)--(0,0)--cycle;
\draw (0,0)--(6,0)--(1.75,4.165)--cycle; \node at (0.6,-0.3)
{$S=(0,0)$}; \node at (6.9,-0.3) {$R=(s+x,0)$}; \node at (3.9,3)
{$D=(p,q)$}; \node at (5,-0.3) {$(s,0)$}; \node at
(4.25,2.0825) {$m_S$};
\node at
(0.5,2.0825) {$m_R$};
\node at
(3,0.2) {$m_P$};
\draw[fill=black] (0,0) circle (0.05); \draw[fill=black] (3,0)
circle (0.05); \draw[fill=black] (0.875,2.0825) circle (0.05);
\draw[fill=black] (5,0) circle (0.05);

\draw[fill=black] (3, 2.95) circle (0.05); \draw[fill=black]
(3.875,2.0825) circle (0.05);
\coordinate
(S) at (0,0);
\coordinate
(R) at (6,0); \coordinate[label=left:$P$] (P) at (1.75,4.165); \node
at (5.6,0.2) {$\delta$}; \draw pic[draw=black, angle
eccentricity=1.2, angle radius=0.7cm]
    {angle=P--R--S};
\draw (1.75,4.165)--(1.75,0); \node at (1.5,2) {$m$};

\draw[fill=gray!20] (8.5,0)--(11.7,0)--(12,2.48)--cycle ;
\draw[dashed] (8.5,0)--(12,2.48)--(15.5,0)--cycle;
\draw (8.5,0)--(14.5,0)--(10.25,4.165)--cycle; \node at (9.1,-0.3)
{$S=(0,0)$}; \node at (15.4,-0.3) {$R=(s+x,0)$}; \node at (11.8,-0.3)
{$(s,0)$}; \node at (12.9,2.58) {$D=(p,q)$};
\draw[fill=black] (8.5,0) circle (0.05); \draw[fill=black] (11.5,0)
circle (0.05); \node at
(12.75,2.0825) {$m_S$};
\node at
(9,2.0825) {$m_R$};
\node at
(11.4,0.2) {$m_P$};

\draw[fill=black] (9.375,2.0825) circle (0.05); \draw[fill=black]
(11.7,0) circle (0.05); \draw[fill=black] (12,2.48) circle (0.05);
\draw[fill=black] (12.375,2.0825) circle (0.05); \coordinate (S) at
(8.5,0); \coordinate  (R) at (14.5,0); \coordinate [label=left:$P$] (P)
at (10.25,4.165); \node at (14.1,0.2) {$\delta$}; \draw
pic[draw=black, angle eccentricity=1.2, angle radius=0.7cm]
    {angle=P--R--S};

\draw (10.25,4.165)--(10.25,0); \node at (10,2) {$m$};

\end{tikzpicture}
\caption{Possible realizations of Case (4) when $\delta<60 ^{\circ}$.}
    \label{fig4.44}
\end{center}
\end{figure}
\noindent
By simple calculation,
\begin{equation}\label{eqP}P=(s+x,0)+\frac{(s+x)}{\sqrt{(p-(s+x))^2+q^2}}(p-(s+x),q).\end{equation}
Denote by $m$ the length of the altitude of $SPR$ belonging to the side $SR$. Then, $m$ is equal to the second coordinate of $P$ (see Figure \ref{fig4.44}). We have
$$m=\frac{q(s+x)}{\sqrt{(p-(s+x))^2+q^2}}.$$

Let us compute the derivative of the function $$f(x)=2t(SPR)=
q(s+x)^2\cdot\big((p-(s+x))^2+q^2\big)^{-\frac{1}{2}}.$$
We obtain
\begin{equation*}
\begin{split}
f'(x)&= 2q(s+x)\big((p-(s+x))^2+q^2\big)^{-\frac{1}{2}}+q(s+x)^2(p-(s+x))\big((p-(s+x))^2+q^2\big)^{-\frac{3}{2}}\\
&=q(s+x)\big((p-(s+x))^2+q^2\big)^{-\frac{3}{2}}\big[(p-(s+x))(2p-(s+x))+2q^2\big]\\
&=
\underbrace{q\left(s+x\right)\left((p-(s+x))^2+q^2\right)^{-\frac{3}{2}}}_{>0}\left[\left(\frac{3}{2}p-(s+x)\right)^2-\frac{1}{4}p^2+2q^2\right].
\end{split}
\end{equation*}

{\em Case (4/a1)}: $q \ge \frac{1}{2}p$. Then, $-\frac{1}{4}p^2+2q^2>0$, and
hence, $f'(x)>0$ for all $x\ge 0$. Thus, $f$ is strictly increasing and since $x$ cannot be negative, $f$
takes its minimum at $x=0$. This means that the area of a special container
of the first kind where $x=0$  (see the triangle with dashed sides on Figure \ref{fig4.44} (4/a)) is smaller than the area of $SPR$ for $x>0$.

{\em Case (4/a2)}: $2p<s+x$. Then $\frac{1}{2}p<(s+x)-\frac{3}{2}p$, so that $(\frac{3}{2}p-(s+x))^2-(\frac{1}{2}p)^2>0$. Again, we have $f'(x)>0$ for all $x\ge 0$ and, as above, we obtain a special container of the first kind whose area is smaller than the area of $SPR$ (see Figure \ref{fig4.44}
(4/a)).

{\em Case (4/b)}: $q<\frac{1}{2}p$ and $2p\ge s+x$. Let $\Delta$ denote the triangle with vertices $(0,0)$, $(p,q)$, and $(2p,0)$. It follows from the inequality $2p\ge s+x$ that $\Delta$ is an isosceles container of the second kind associated with $ABC$ (see Figure \ref{fig4.44}
(4/b)). We show that $\Delta$ has smaller area than $SPR$.
To prove this, we have to verify that
$$t(\Delta)=pq<q\frac{(s+x)^2}{2\sqrt{(p-(s+x))^2+q^2}}=t(SPR).$$
Using our assumption that $(p,q)\in m_SP$, we obtain $$p\le
s+x+\frac{(s+x)(p-(s+x))}{2\sqrt{(p-(s+x))^2+q^2}}.$$ The right-hand side of the last inequality is the first coordinate of the midpoint $m_S$ of $PR$, where $P$ is given by formula \eqref{eqP} and $R=(s+x,0)$.  Thus, it is
sufficient to show that
$$(s+x)+\frac{(s+x)(p-(s+x))}{2\sqrt{(p-(s+x))^2+q^2}}<\frac{(s+x)^2}{2\sqrt{(p-(s+x))^2+q^2}},$$
which reduces to $3p^2+4q^2<4(s+x)p$. The last inequality holds, because $p<s+x$ and, by our assumption, $2q\le p$. This completes the proof of Theorem \ref{t0.2}. %and, hence, of Theorem
%\ref{t0.1}.
\qed
\medskip

\noindent{\bf Proof of Theorem \ref{prop1.3}.} By Theorem \ref{t0.2}, every
minimum area isosceles container for a triangle $ABC$ is one of its
special containers. By Lemma \ref{p2}, it must be a special container of the
first or the second kind.

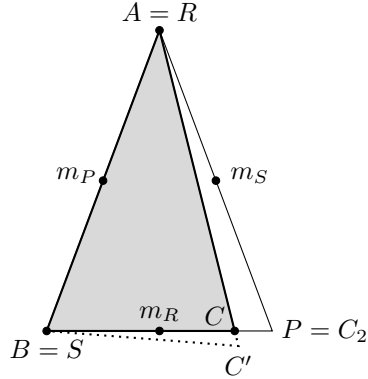
\begin{figure}[!ht]
\begin{center}
    \begin{tikzpicture}
\draw (0,-11) node[anchor=north]{$B=S$}
  -- (3,-11) node[anchor=west]{$P=C_2$}
  -- (1.5,-7) node[anchor=south]{$A=R$}
  -- cycle;

 \draw[fill=gray!30] (0,-11)-- (1.5,-7) --(2.5,-11)--cycle;
%\draw (0,-7.1) node[anchor=south]{$(10)$};
 \draw[fill=black] (1.5,-11) circle (0.05);
 \draw[fill=black] (2.25,-9) circle (0.05);
 \draw[fill=black] (0.75,-9) circle (0.05);

\draw (1.5,-11) node[anchor=south]{$m_{R}$}; \draw (2.7,-8.7)
node[anchor=north]{$m_{S}$}; \draw (0.4,-8.7)
node[anchor=north]{$m_{P}$};

\draw[fill=black] (0,-11) circle (0.05);
 \draw[fill=black] (1.5,-7) circle (0.05);
 \draw[fill=black] (2.5,-11) circle (0.05);
\draw[dotted, thick] (2.5,-11)--(2.55,-11.2)--(0,-11);

\draw (2.5,-10.8) node[anchor=east]{$C$}; \draw (2.55,-11.2)
node[anchor=north]{$C'$};

\draw[thick] (0,-11)-- (1.5,-7) -- (2.5,-11)-- cycle;

    \end{tikzpicture}
\end{center}
    \caption{An acute minimum area isosceles container cannot be a special container of the second kind.}
    \label{figfirst}
\end{figure}

Suppose for a contradiction that a minimum area special container $SPR$ associated
with the triangle $ABC$ is acute, but it is of the second kind. Assume without
loss of generality that $S=B$ and $R=A$. (The other cases can be treated in a similar manner.) By our notation, $P=C_2$. We prove that $t(ABC')<t(ABC_2)=t(SPR)$ (see Figure \ref{figfirst}).

Indeed, we have $|AB|=|AC_2|=|AC'|$ and $\sa C'AB<\sa C_2AB$. Since
$SPR=ABC_2$ is acute, it follows that $t(ABC')<t(ABC_2)$. \qed
\medskip

\begin{cor}\label{cor2}
If a minimum area isosceles container for $ABC$ is a special container of the second kind, then it must be $AB_1C$.
\end{cor}

\noindent{\bf Proof.}
By Theorem \ref{prop1.3}, if a special container of the second kind has minimum area, then it has to be non-acute. If $ABC_2$ is non-acute, then $ABC_1$ is obtuse and $t(ABC_2)>t(ABC_1)$, because $|AB|=|AC_1|=|AC_2|$ and $\sa ABC_1>\sa BAC_2\ge90 ^{\circ}$. On the other hand, as $AB_1C$ and $ABC_1$ share a base angle at $A$ and $b<c$, it follows that $t(ABC_1)>t(AB_1C)$.   \hfill    $\Box$

\section{Quantitative results---Proofs of Theorems~\ref{peq3} and~\ref{thmrarea} (\rm{a})}\label{uj}

\noindent{\bf Proof of Theorem \ref{peq3}.} By Theorem \ref{t0.2}, a minimal area isosceles container for $ABC$ is a special container associated with $ABC$. In view of Lemma~\ref{p2}, it must be a special container of the first or second kind. By Lemma
\ref{p0.1} (a) and Corollary \ref{cor2}, among special containers of the first kind, it is
enough to consider $ABC'$ and $AB'C$, and among special containers of the second kind, only  $AB_1C$. These immediately show that every triangle $ABC$ admits at most 3 minimum area isosceles containers.

If $ABC$ is an obtuse or right triangle,
then $t(AB'C) >t(AB_1C)$. Indeed, in this case $|AC|=|CB'|=|CB_1|$, both $AB'C$ and $AB_1C$ are obtuse or right triangles, and their apex angles satisfy $\sa ACB' <\sa ACB_1$. Thus, there are only two candidates for a minimum area isosceles container: $ABC'$ and $AB_1C$.

If $ABC$ is an acute triangle
and it has $3$ minimum area isosceles containers, then
$t(AB'C)=t(ABC')=t(AB_1C)$.
Since $t(BCB_1)=t(BCC')$, we obtain
\begin{equation}\label{eq1} (c-b) \sin(\alpha+\beta)=b\sin(\beta - \alpha).
\end{equation}
Note that this equation also holds when $ABC$ is obtuse.

It follows from equation \eqref{eq1} that $\frac{c}{b}=\frac{2
\sin(\beta)\cos(\alpha)}{\sin(\alpha+\beta)}$. By Lemma \ref{p0.1} (b), the equation
$t(AB'C)=t(ABC')$ reduces to $\frac{c}{b}=\frac{b}{a}$. Thus,
$\frac{\sin(\beta)}{\sin(\alpha)}=\frac{2
\sin(\beta)\cos(\alpha)}{\sin(\alpha+\beta)}$, so that
$\sin(\alpha+\beta)=\sin(2\alpha)$. Therefore, either $\alpha=\beta$,
which is impossible, or $180 ^{\circ}-(\alpha+\beta)=\gamma=2\alpha$.

It follows from $\frac{c}{b}=\frac{b}{a}$ that
\begin{equation}\label{eq3}
    \frac{\sin(2\alpha)}{\sin(3\alpha)}=\frac{\sin(3\alpha)}{\sin(\alpha)}.
\end{equation}
Since $ABC$ is acute, its smallest angle is $\alpha$, and $180 ^{\circ}-3\alpha=\beta<\gamma=2\alpha$, we
have that $36 ^{\circ}<\alpha<45 ^{\circ}$. Simple analysis shows that equation
\eqref{eq3} has exactly one solution $\alpha^*$ in the interval $[36 ^{\circ},45 ^{\circ}]$. It can be approximated by computer. The other two angles of the corresponding triangle are $\beta^*=180 ^{\circ}-3\alpha^*$ and $\gamma^*=2\alpha^*$.  \qed

\medskip
\begin{ex}\label{ex1}
By the proof of Theorem \ref{peq3}, any minimal area isosceles container for $ABC$ is either $AB'C$, or $ABC'$, or $AB_1C$.
Here, we construct a family of acute triangles $ABC$ whose only minimal area isosceles containers are special containers of the second kind, i.e., $AB_1C$. Moreover, $AB_1C$ is obtuse.

Let $\alpha>\alpha^*$ and $90 ^{\circ}>\gamma>2\alpha>\gamma^*$. Then, denoting by $\beta=180^{\circ}-\alpha-\gamma$ we obtain that
\[\frac{\sin(\gamma)}{\sin(\beta)}>\frac{\sin(\gamma^*)}{\sin(\beta^*)}=
\frac{\sin(\beta^*)}{\sin(\alpha^*)}>
\frac{\sin(\beta)}{\sin(\alpha)},\] which implies that $t(AB'C)<t(ABC')$.
The triangles $AB_1C$ and $AB'C$ are isosceles with legs of length $b$,
so it is enough to show that $\sin(\sa ACB_1)=\sin(180 ^{\circ}-2 \alpha)<
\sin( \gamma )=\sin(\sa ACB')$. However, this follows from the inequalities
$90 ^{\circ}>\gamma>2\alpha$. The base angle of $AB_1C$ satisfies $\alpha<45 ^{\circ}$, so that
$AB_1C$ is obtuse.

\end{ex}

\noindent{\bf Proof of Theorem \ref{thmrarea} (a):}
Let $r^*$ denote the supremum of the ratio of the area of a minimum area isosceles container of a triangle to the area of the triangle itself. In view of Theorem \ref{peq3}, we have $$r^*=\sup_{\textrm{triangle } ABC}\min\Big(
\frac{t(AB'C)}{t(ABC)},\frac{t(ABC')}{t(ABC)}, \frac{t(AB_1C)}{t(ABC)}\Big).$$

If $\beta\ge 45 ^{\circ}$, then $\sin(\beta)\ge\frac{1}{\sqrt{2}}$. Using the proof of Lemma \ref{p0.1}(b) and
the law of sines, we obtain $$\frac{t(ABC')}{t(ABC)}=
\frac{c}{b}=\frac{\sin(\gamma)}{\sin(\beta)}\le
\frac{1}{\frac{1}{\sqrt{2}}}=\sqrt{2}.$$ Equality holds here if and only if $\beta=
45 ^{\circ}$ and $\gamma= 90 ^{\circ}$, in which case $ABC$ is an isosceles triangle and
the ratio of the area of the minimum isosceles container to the area of $ABC$ is $1$.

If $\beta < 45 ^{\circ}$, then $\gamma>90 ^{\circ}$. Hence, $ABC$ is obtuse and, by the proof of Theorem \ref{peq3}, the minimum
area isosceles container is $ABC'$ or $AB_1C$. For fixed
$\beta$ and $c$, we can express the ratios of the areas as
functions of $\alpha$. Let $$f(\alpha)=\frac{t(ABC')}{t(ABC)}=
\frac{c}{b}$$ and
$$g(\alpha)=\frac{t(AB_1C)}{t(ABC)}=\frac{2b\cos(\alpha)}{b\cos(\alpha)+a\cos(\beta)}=\frac{1}{\frac{1}{2}+\frac{\tan(\alpha)}{2\tan(\beta)}},$$
where $0<\alpha<\beta$. Obviously, $f(\alpha)$ is strictly
increasing and $g(\alpha)$ is strictly decreasing on the open interval $(0,\beta)$, and
both functions are continuous. We have $$\lim_{\alpha\to 0+} f(\alpha)=1, \ \
1<\lim_{\alpha\to \beta-}f(\alpha)<2,$$ $$\lim_{\alpha\to
0+} g(\alpha)=2, \ \ \lim_{\alpha\to \beta-}g(\alpha)=1.$$
Therefore, the graphs of $f$ and $g$ intersect at a unique point $z$. Thus,
$\max\limits_{0 ^{\circ}<\alpha<\beta}^{}({\min\left(f(\alpha),g(\alpha)\right)})= f(z)=g(z)$, which implies $t(ABC')=t(AB_1C)$. This means that $t(BCB_1)=t(BCC')$, so that equation \eqref{eq1} above holds. Using the law of sines, we obtain $(\frac{c}{b})^2=2
\cos(z)< 2$ and, hence, $\frac{c}{b}< \sqrt{2}$.
If $\beta\to 0$, then $z\to 0$ and $c/b \to \sqrt{2}$.
This implies that $r^*=\sqrt{2}$, but the supremum is not realized by any triangle $ABC$.
\qed

\bigskip

\noindent{\bf Acknowledgement.}
The authors would like to express their gratitude to the anonymous referee for her useful suggestions. Research by the first author was supported by a Premium Postdoctoral Fellowship of the Hungarian Academy of Science, and by NKFIH (Hungarian National Research, Development
and Innovation Office) grant K-124749. Research by the second author was partially supported by NKFIH grants K-176529 and KKP-133864, Austrian Science Fund grant Z 342-N31, and by the Ministry of Education and Science of the Russian Federation in the framework of MegaGrant No.\ 075-15-2019-1926. Research by the third author was partially supported by NKFIH grant K-115799.

\medskip

Gergely Kiss,
Alfr\'ed R\'enyi Institute of Mathematics, POB 127, Budapest H-1364, Hungary \\
{\it E-mail address:} kigergo57@gmail.com
\\

J\'anos Pach,
Alfr\'ed R\'enyi Institute of Mathematics, POB 127, Budapest H-1364, Hungary and MIPT, Dolgoprudny, Moscow Oblast 141701, Russia \\
{\it E-mail address:} pachjanos@gmail.com
\\

G\'abor Somlai, L\'or\'and E\"otv\"os University, P\'azm\'any P\'eter
s\'et\'any 1/C,
Budapest H-1117, Hungary\\
{\it E-mail address:} zsomlei@gmail.com
\\
\end{document}